\documentclass[leqno,centertags]{amsart}
\usepackage{amsfonts}
\usepackage{amssymb}
\usepackage{graphicx}
\usepackage{amscd}
\usepackage{endnotes}

\theoremstyle{plain}

\numberwithin{equation}{section}

\input{tcilatex}

\begin{document}

\begin{center}
{\LARGE Existence and higher arity iteration for total asymptotically
nonexpansive mappings in uniformly convex hyperbolic spaces\bigskip }

Hafiz Fukhar-ud-din$^{a,b},$ Amna Kalsoom$^{b}$ and Muhammad Aqeel Ahmad Khan%
$^{b,\ast }$

$^{a}$Department of Mathematics and Statistics, King Fahd University of
Petroleum and Minerals, Dhahran, 3126, Saudi Arabia

$^{b}$Department of Mathematics, The Islamia University of Bahawalpur,

Bahawalpur, 63100, Pakistan

December 19, 2013
\end{center}

\textbf{Abstract}: This paper provides a fixed point theorem and iterative
construction of a common fixed point for a general class of nonlinear
mappings in the setup of uniformly convex hyperbolic spaces. We translate a
multi-step iteration, essentially due to Chidume and Ofoedu \cite{Chidume
JMAA}\ in such setting for the approximation of common fixed points of a
finite family of total asymptotically nonexpansive mappings. As a
consequence, we establish strong and $\triangle $-convergence results which
extend and generalize various corresponding results announced in the current
literature.\newline

\section{Introduction}

Let $(X,d)$ be a metric space. A geodesic from $x$ to $y$ in $X$ is a
mapping $\theta :[0,I]\rightarrow 
\mathbb{R}
$ such that%
\begin{equation*}
d(\theta (t),\theta (s)=\left\vert t-s\right\vert ,~\forall ~t,s\in \lbrack
0,I].
\end{equation*}

The above characteristics is referred as a constant speed of $\theta $, the
parametrization of $\theta $ w.r.t. the arc length or distance preservation
of $\theta $. A geodesic or metric segment is the image $\Theta $ of a
geodesic $\theta :[0,I]\rightarrow 
\mathbb{R}
$ in $X.$ The points $x=\theta (0)$ and $y=\theta (I)$ are called the end
points or the extreme (maximal or minimal) lengths of that segment. A metric
space is coined as a unique geodesic space if every two points of the metric
space $(X,d)$ are joint by a unique geodesic segment.

A geodesic is an arbitrary curve between two points on a surface. It plays
an important role in real world applications, for example geometric
designing of various shapes and differential geometry. The class of
hyperbolic spaces introduced by Kohlenbach \cite{Kohlenbach} is an important
example of (uniquely-) geodesic space as well as prominent among various
other notions of hyperbolic spaces in the current literature, see for
example \cite{Goebel Kirk,Goebel Reich,Reich Shafrir}. The study of
hyperbolic spaces has been largely motivated and dominated by questions
about hyperbolic groups, one of the main object of study in geometric group
theory. We remark that non-positively curved spaces, such as hyperbolic
spaces, play a significant role in many branches of applied mathematics.

Throughout this paper, we work in the setting of hyperbolic spaces
introduced by Kohlenbach \cite{Kohlenbach}:

A hyperbolic space \cite{Kohlenbach} is a metric space $(X,d)$ together with
a mapping $W:$ $X^{2}\times \lbrack 0,1]\rightarrow X$ satisfying%
\begin{equation*}
\begin{array}{l}
\text{(W1) }d(u,W(x,y,\alpha ))\leq \alpha d(u,x)+(1-\alpha )d(u,y) \\ 
\text{(W2) }d(W(x,y,\alpha ),W(x,y,\beta ))=\left\vert \alpha -\beta
\right\vert d(x,y) \\ 
\text{(W3) }W(x,y,\alpha )=W(y,x,(1-\alpha )) \\ 
\text{(W4) }d(W(x,z,\alpha ),W(y,w,\alpha ))\leq (1-\alpha )d(x,y)+\alpha
d(z,w)\newline
\end{array}%
\end{equation*}%
for all $x,y,z,w\in X$ and $\alpha ,\beta \in \lbrack 0,1].$

The class of hyperbolic spaces in the sense of Kohlenbach \cite{Kohlenbach}
contains all normed linear spaces and convex subsets thereof as well as
Hadamard manifolds and $CAT(0)$ spaces in the sense of Gromov. An important
example of a hyperbolic space is the open unit ball $B_{H}$ in a real
Hilbert space $H$ is as follows:

let $B_{H}$ be the open unit ball in $H.$ Then%
\begin{equation*}
k_{B_{H}}\left( x,y\right) :=\arg \tanh \left( 1-\sigma (x,y)\right) ^{\frac{%
1}{2}},
\end{equation*}%
where%
\begin{equation*}
\sigma (x,y)=\frac{\left( 1-\left\Vert x\right\Vert ^{2}\right) \left(
1-\left\Vert x\right\Vert ^{2}\right) }{\left\vert 1-\left\langle
x,y\right\rangle \right\vert ^{2}}\text{ for all }x,y\in B_{H},
\end{equation*}%
defines a metric on $B_{H}~$(also known as Kobayashi distance). The open
unit ball $B_{H}$ together with this metric is coined as Hilbert ball. Since 
$(B_{H},k_{B_{H}})$ is a unique geodesic space, so one can define a
convexity mapping $W$ for the corresponding hyperbolic space $%
(B_{H},k_{B_{H}},W).$ This space is of significant importance for fixed
point theory of holomorphic mappings as the said class of mappings is $%
k_{B_{H}}$-nonexpansive in $(B_{H},k_{B_{H}},W).$ A metric space $(X,d)$
satisfying only (W1) is a convex metric space introduced by Takahashi \cite%
{Takahashi}. A subset $K$ of a hyperbolic space $X$ is convex if $%
W(x,y,\alpha )\in K$ for all $x,y\in K$ and $\alpha \in \lbrack 0,1].$ For
more on hyperbolic spaces and a comparison between different notions of
hyperbolic spaces present in the literature, we refer to \cite[p.384]%
{Kohlenbach}.

It is worth to mention that fixed point theory of nonexpansive mappings
(i.e., $d(Tx,Ty)\leq d(x,y),~$for $x,y\in K$) and its various
generalizations majorly depends on the geometrical characteristics of the
underlying space. The class of nonexpansive mappings enjoys the fixed point
property (FPP) and the approximate fixed point property (AFPP) in various
setting of spaces. Moreover, it is natural to extend such powerful results
to generalized nonexpansive mappings as a means of testing the limit of the
theory of nonexpansive mappings. FPP and even AFPP, in a nonlinear domain,
of various generalizations of nonexpansive mappings are still developing.
The class of hyperbolic spaces is endowed with rich geometric structures for
different results with applications in topology, graph theory, multivalued
analysis and metric fixed point theory. An important ingredient for metric
fixed point theory of nonexpansive mappings is \textit{uniform convexity}.

A hyperbolic space $(X,d,W)$ is uniformly convex \cite{L-JMAA} if for all $%
r>0$ and $\epsilon \in (0,2]$ there exists a $\delta \in (0,1]$ such that
for all $u,x,y\in X,$ we have%
\begin{equation*}
d\left( W(x,y,\frac{1}{2}),u\right) \leq (1-\delta )r
\end{equation*}%
provided $d(x,u)\leq r,d(y,u)\leq r$ and $d(x,y)\geq \epsilon r.$

A mapping $\eta :(0,\infty )\times (0,2]\rightarrow (0,1],$ which provides
such a $\delta =\eta (r,\epsilon )$ for given $r>0$ and $\epsilon \in (0,2],$%
\ is known as a modulus of uniform convexity of $X$ . We call $\eta $
monotone if it decreases with $r$ (for a fixed $\epsilon $), $i.e.,~\forall
\epsilon >0,\forall r_{2}\geq r_{1}>0~(\eta \left( r_{2},\epsilon \right)
\leq \eta \left( r_{1},\epsilon \right) ).$ $CAT(0)$ spaces are uniformly
convex hyperbolic spaces with modulus of uniform convexity $\eta (r,\epsilon
)=\frac{\epsilon ^{2}}{8}$ \cite{L-JMAA}. Therefore , the class of uniformly
convex hyperbolic spaces includes both uniformly convex normed spaces and $%
CAT(0)$ spaces as special cases.

Metric fixed point theory of nonlinear mappings in a general setup of
hyperbolic spaces is a fascinating field of research in nonlinear functional
analysis. Moreover, iteration schemas are the only main tool to study fixed
point problems of nonexpansive mappings and its various generalizations. In
2006, Alber et al. \cite{Alber FPTA} introduced a unified and generalized
notion of a class of nonlinear mappings in Banach spaces, which can be
introduced in the general setup of hyperbolic spaces as follows:

Let $T:K\rightarrow K$ be a self-mapping, then $T$ is called total
asymptotically nonexpansive mappings if there exists nonnegative real
sequences $\{k_{n}\}$ and $\{\varphi _{n}\}$ with $k_{n}\overset{%
n\rightarrow \infty }{\rightarrow }0,\varphi _{n}\overset{n\rightarrow
\infty }{\rightarrow }0$ and a strictly increasing continuous function $\xi :%
\mathbb{R}
^{+}\rightarrow 
\mathbb{R}
^{+}$ with $\xi (0)=0,$ such that 
\begin{equation*}
d(T^{n}x,T^{n}y)\leq d(x,y)+k_{n}\xi \left( d(x,y)\right) +\varphi _{n}\text{
\ for all \ }x,y\in K,~n\geq 1.
\end{equation*}%
\textbf{Example 1.1.} (i) Let $X:=%
\mathbb{R}
$ and $K=[0,\infty ),$ then the mapping $T:K\rightarrow K$ be defined by $%
Tx=\sin x$ is total asymptotically nonexpansive.\newline
(ii) Let $K=\left[ -\frac{1}{\pi },\frac{1}{\pi }\right] ,$ then the mapping 
$T:K\rightarrow K$ be defined by $Tx=kx\sin \frac{1}{x},$ where $k\in (0,1)$
is total asymptotically nonexpansive.\newline
(iii) Let $K=\{x:=(x_{1},x_{2},\cdots ,x_{n},\cdots )~|~x_{1}\leq 0,x_{i}\in 
\mathbb{R}
$ $i\geq 2\}[0,\infty )$ be a nonempty subset of $X=l^{2}$ with the norm $%
\left\Vert \cdot \right\Vert $ defined as:%
\begin{equation*}
\left\Vert x\right\Vert =\sqrt{\tsum_{i=1}^{\infty }x_{i}^{2}}.
\end{equation*}%
Let $T:K\rightarrow K$ be defined by%
\begin{equation*}
T(x)=(0,4x_{2},0,0,0,\cdots ),
\end{equation*}%
then $T$ is total asymptotically nonexpansive.\newline
(iv) Let $X=l^{2}$and $K$ be a unit ball in $l^{2}$.\ Let $T:K\rightarrow K$
be a mapping defined by%
\begin{equation*}
T(x)\rightarrow (0,x_{1}^{2},a_{2}x_{2},\cdots ,a_{n}x_{n},\cdots ),\text{
for all }x\in l^{2}
\end{equation*}%
where $x:=(x_{1},x_{2},\cdots ,x_{n},\cdots )$ and $\left\{ a_{i}\right\} $
is a sequence in $(0,1)$ such that $\tprod_{i=2}^{\infty }a_{i}=2.$ Then $T$
is total asymptotically nonexpansive.

The class of total asymptotically nonexpansive mappings and asymptotically
nonexpansive mappings have been studied extensively in the literature, see
for example, see for example \cite{Chang AMC,Fukhar Khan CAMWA,FKK 2012,Khan
JIA,Khan Fukhar FPTA,Yildirim-Khan} and the references cited therein. It is
worth mentioning that the results established for total asymptotically
nonexpansive mappings are applicable to the mappings associated with the
class of asymptotically nonexpansive mappings and which are extensions of
nonexpansive mappings.

It is remarked that the iterative construction of common fixed points of a
finite family of asymptotically quasi-nonexpansive mappings in a Banach
space through higher arity of an iteration is essentially due to Khan et al. 
\cite{KDF JMAA}(see also \cite{Kuhfittig} for the case of nonexpansive
mappings). This iteration was further generalized by Khan and Ahmad \cite%
{Khan Ahmad 2010} and Khan et al. \cite{KKF NA} to the setup of convex
metric spaces and $CAT(0)$ spaces, respectively. Moreover, Chidume and
Ofoedu \cite{Chidume JMAA} introduced a general iteration schema for a
finite family of total asymptotically nonexpansive mappings in Banach
spaces, which we adopt in the setting of hyperbolic spaces as follows:%
\begin{equation}
\begin{array}{l}
x_{1}\in K, \\ 
x_{n+1}=W\left( x_{n},T_{1}^{n}x_{n},\alpha _{n}\right) \text{, if }%
m=1,n\geq 1, \\ 
x_{1}\in K, \\ 
x_{n+1}=W\left( x_{n},T_{1}^{n}y_{1n},\alpha _{n}\right) , \\ 
y_{1n}=W\left( x_{n},T_{2}^{n}y_{2n},\alpha _{n}\right) , \\ 
y_{2n}=W\left( x_{n},T_{3}^{n}y_{3n},\alpha _{n}\right) , \\ 
\multicolumn{1}{c}{\vdots} \\ 
y_{(m-2)n}=W\left( x_{n},T_{m-1}^{n}y_{(m-1)n},\alpha _{n}\right) , \\ 
y_{(m-1)n}=W\left( x_{n},T_{m}^{n}x_{n},\alpha _{n}\right) \text{ , if }%
m\geq 2,n\geq 1,%
\end{array}
\tag{1.1}
\end{equation}%
where $\left\{ \alpha _{n}\right\} _{n=1}^{\infty }$ is a sequence in $[0,1]$
bounded away from $0$ and $1.$

The purpose of this paper is to establish fixed point result along with the
iterative construction of common fixed point of a finite family of mappings
in uniformly convex hyperbolic spaces. We therefore, establish results
concerning strong convergence and $\triangle $-convergence results of
iteration (1.1). Our convergence results can be viewed not only as analogue
of various existing results but also improve and generalize various results
in the current literature.

\section{Preliminaries and Some Auxiliary Lemmas}

We start this section with notion of asymptotic center -- essentially due to
Edelstein \cite{Edelstein BAMS}-- of a sequence which is not only useful in
proving that a fixed point exists but also plays a key role to define the
concept of $\triangle $-convergence in hyperbolic spaces. In 1976, Lim \cite%
{Lim} introduced the concept of $\triangle $-convergence in the general
setting of metric spaces. In 2008, Kirk and Panyanak \cite{Kirk Panyanak}
further analyzed this concept in geodesic spaces. They showed that many
Banach space results involving weak convergence have precise analogue
version of $\triangle $-convergence in geodesic spaces.

Let $\{x_{n}\}$ be a bounded sequence in a hyperbolic space $X$. For $x\in X$%
, define a continuous functional $r(.,\{x_{n}\}):X\rightarrow \lbrack
0,\infty )$ by:%
\begin{equation*}
r(x,\{x_{n}\})=\limsup_{n\rightarrow \infty }d(x,x_{n}).
\end{equation*}%
The asymptotic radius and asymptotic center of the bounded sequence $%
\{x_{n}\}$ with respect to a subset $K$ of $X$ is defined and denoted as:%
\begin{equation*}
r_{K}(\{x_{n}\})=\inf \{r(x,\{x_{n}\}):x\in X\},
\end{equation*}%
and%
\begin{equation*}
A_{K}(\{x_{n}\})=\{x\in K:r(x,\{x_{n}\})\leq r(y,\{x_{n}\})\text{ for all }%
y\in K\},
\end{equation*}%
respectively.\newline
\textbf{Lemma 2.1} (\cite{L-Proceeding}). \textit{Let }$(X,d,W)$\textit{\ be
a complete uniformly convex hyperbolic space with monotone modulus of
uniform convexity }$\eta $\textit{. Then every bounded sequence }$\{x_{n}\}$%
\textit{\ in }$X$\textit{\ has a unique asymptotic center with respect to
any nonempty closed convex subset }$K$\textit{\ of }$X$\textit{.}

Recall that a sequence $\{x_{n}\}$ in $X$ is said to $\triangle -$converge
to $x\in X$ if $x$ is the unique asymptotic center of $\{u_{n}\}$ for every
subsequence $\{u_{n}\}$ of $\{x_{n}\}.$ In this case, we write $\triangle
-\lim_{n}x_{n}=x$ and call $x$ as $\triangle -$limit of $\{x_{n}\}.$ A
mapping $T:K\rightarrow K$ is \textit{semi-compact} if every bounded
sequence $\{x_{n}\}\subset K$ satisfying $d(x_{n},Tx_{n})\rightarrow 0,$ has
a convergent subsequence.

We now list some useful lemmas as well as establish some auxiliary results
required in the sequel.\newline
\textbf{Proposition 2.2} (\cite{Kohlenbach JEMS}). \textit{Let }$(X,d,W)$%
\textit{\ be a complete uniformly convex hyperbolic space with monotone
modulus of uniform convexity }$\eta $\textit{. The intersection of any
decreasing sequence of nonempty bounded closed convex subsets of }$X$\textit{%
\ is nonempty.}\newline
\textbf{Lemma 2.3.}(\cite{L-JMAA}) \textit{Let }$(X,d,W)$\textit{\ be a
uniformly convex hyperbolic space with monotone modulus of uniform convexity 
}$\eta .$\textit{\ For }$r>0,~\epsilon \in (0,2],~a,x,y\in X$ and $\lambda
\in \left[ 0,1\right] ,$\textit{\ the inequalities}%
\begin{equation*}
d(x,a)\leq r,~d(y,a)\leq r\text{ and }d(x,y)\geq \epsilon r
\end{equation*}%
\textit{imply}%
\begin{equation*}
d\left( W(x,y,\lambda ),a\right) \leq \left( 1-2\lambda \left( 1-\lambda
\right) \eta \left( r,\epsilon \right) \right) r.
\end{equation*}%
\textbf{Lemma 2.4(}\cite{Chang}). \textit{Let }$\left\{ a_{n}\right\}
,\left\{ b_{n}\right\} $\textit{\ and }$\left\{ c_{n}\right\} $\textit{\ be
sequences of non-negative real numbers such that}$~\sum_{n=1}^{\infty
}b_{n}<\infty $ and $\sum_{n=1}^{\infty }c_{n}<\infty .$\textit{\ If~}$%
a_{n+1}\ \leq (1+b_{n})a_{n}+c_{n},$\textit{\ }$n\geq 1,~$\textit{then }$%
\lim_{n\rightarrow \infty }a_{n}$\textit{\ exists.}\newline
\textbf{Lemma 2.5(}\cite{KFK FPTA}\textbf{)}. \textit{Let }$(X,d,W)$\textit{%
\ be a uniformly convex hyperbolic space with monotone modulus of uniform
convexity }$\eta .$\textit{~Let }$x\in X$\textit{\ and }$\{\alpha _{n}\}$%
\textit{\ be a sequence in }$[a,b]$ for some $a,b\in (0,1).$ \textit{If }$%
\{x_{n}\}$ and $\{y_{n}\}$\textit{\ are sequences in }$X$\textit{\ such that 
}$\limsup_{n\longrightarrow \infty }d(x_{n},x)\leq
c,~\limsup_{n\longrightarrow \infty }d(y_{n},x)\leq c$\textit{\ and }$%
\lim_{n\longrightarrow \infty }d(W(x_{n},y_{n},\alpha _{n}),x)=c$\textit{\
for some }$c\geq 0,$\textit{\ then }$\lim_{n\rightarrow \infty
}d(x_{n},y_{n})=0.$\newline
\textbf{Lemma 2.6(}\cite{KFK FPTA}\textbf{)}. \textit{Let }$K$\textit{\ be a
nonempty closed convex subset of a uniformly convex hyperbolic space and }$%
\{x_{n}\}$\textit{\ a bounded sequence in }$K$\textit{\ such that }$%
A_{K}(\{x_{n}\})=\{y\}$ and $r_{K}(\{x_{n}\})=\rho $\textit{. If }$\{y_{m}\}$%
\textit{\ is another sequence in }$K$\textit{\ such that }$%
\lim_{m\rightarrow \infty }r(y_{m},\{x_{n}\})=\rho ,$ \textit{then }$%
\lim_{m\rightarrow \infty }y_{m}=y.$\newline
\textbf{Lemma 2.7.} \textit{Let }$K$\textit{\ be a nonempty closed and
convex subset of a hyperbolic space }$X$\textit{\ and let }$\left\{
T_{i}\right\} _{i=1}^{m}:K\rightarrow K$ \textit{be a finite family of total
asymptotically nonexpansive mappings with sequences }$\{k_{in}\}$ and $%
\{\varphi _{in}\},$ $n\geq 1,\ i=1,2,\cdots ,m$ \textit{such that }$F=\cap
_{i=1}^{m}F\left( T_{i}\right) \neq \emptyset .$ \textit{For} $i=1,2,\cdots
,m$ \textit{if the following conditions are satisfied:}\newline
(C1) $\sum\limits_{n=1}^{\infty }k_{in}<\infty $ \textit{and}$\
\sum\limits_{n=1}^{\infty }\varphi _{in}<\infty ;$\newline
(C2) \textit{there exists constants }$M_{i},\ M_{i}^{\ast }>0$ \textit{such
that }$\zeta _{i}\left( \lambda _{i}\right) \leq M_{i}^{\ast }\lambda _{i}$%
\textit{\ for all }$\lambda _{i}\geq M_{i}$\textit{. Then the sequence }$%
\left\{ x_{n}\right\} $\textit{\ given by }$\left( 1.1\right) $\textit{\ is
bounded and }$\lim\limits_{n\rightarrow \infty }d\left( x_{n},p\right) \ $%
\textit{exists for each }$p\in F$\textit{.}\newline
\textbf{Proof.} Note that `$\xi $' is an increasing function, therefore $\xi
\left( r\right) \leq \xi \left( M\right) $ for $r\leq M.$ Moreover $\xi
\left( r\right) \leq rM_{0}$ for $r\geq M.$ In either case, this implies
that 
\begin{equation*}
\xi \left( r\right) \leq \xi \left( M\right) +rM_{0}.
\end{equation*}%
Let $p\in F$ and $m=1$ in (1.1), we have%
\begin{eqnarray*}
d\left( x_{n+1},p\right) &=&d\left( W\left( x_{n},T_{1}^{n}x_{n},\alpha
_{n}\right) ,p\right) \\
&\leq &\left( 1-\alpha _{n}\right) d\left( x_{n},p\right) +\alpha
_{n}d\left( T_{1}^{n}x_{n},p\right) \\
&\leq &\left( 1-\alpha _{n}\right) d\left( x_{n},p\right) +\alpha
_{n}\left\{ d\left( x_{n},p\right) +k_{1n}\zeta _{1}\left( d\left(
x_{n},p\right) \right) +\varphi _{1n}\right\} \\
&\leq &d\left( x_{n},p\right) +\alpha _{n}k_{1n}\zeta _{1}\left(
M_{1}\right) +\alpha _{n}k_{1n}M_{1}^{\ast }d\left( x_{n},p\right) +\alpha
_{n}\varphi _{1n} \\
&\leq &\left( 1+a_{1}k_{1n}\right) d\left( x_{n},p\right) +a_{1}\left(
k_{1n}+\varphi _{1n}\right) ,
\end{eqnarray*}%
for some constant $a_{1}>0.$\newline
If $m=2$ in (1.1), then we have the following estimates:%
\begin{equation}
\begin{split}
d\left( x_{n+1},p\right) & =d\left( W\left( x_{n},T_{1}^{n}y_{1n},\alpha
_{n}\right) ,p\right) \\
& \leq \left( 1-\alpha _{n}\right) d\left( x_{n},p\right) +\alpha
_{n}d\left( T_{1}^{n}y_{1n},p\right) \\
& \leq \left( 1-\alpha _{n}\right) d\left( x_{n},p\right) +\alpha
_{n}\left\{ d\left( y_{1n},p\right) +k_{1n\ }\zeta _{1}\left( d\left(
y_{1n},p\right) \right) +\varphi _{1n}\right\} \\
& \leq \left( 1-\alpha _{n}\right) d\left( x_{n},p\right) +\alpha
_{n}d\left( y_{1n},p\right) +\alpha _{n}k_{1n}\zeta _{1}\left( M_{1}\right)
\\
& +\alpha _{n}k_{1n\ }M_{1}^{\ast }d\left( y_{1n},p\right) +\alpha
_{n}\varphi _{1n},
\end{split}
\label{2.1}
\end{equation}%
and%
\begin{equation}
\begin{split}
d\left( y_{1n},p\right) & =d\left( W\left( x_{n},T_{2}^{n}x_{n},\alpha
_{n}\right) ,p\right) \\
& \leq \left( 1-\alpha _{n}\right) d\left( x_{n},p\right) +\alpha
_{n}d\left( T_{2}^{n}x_{n},p\right) \\
& \leq \left( 1-\alpha _{n}\right) d\left( x_{n},p\right) +\alpha
_{n}\left\{ d\left( x_{n},p\right) +k_{2n\ }\zeta _{2}\left( d\left(
x_{n},p\right) \right) +\varphi _{2n}\right\} \\
& \leq d\left( x_{n},p\right) +\alpha _{n}k_{2n}\xi _{2}\left( M_{2}\right)
+\alpha _{n}k_{2n\ }M_{2}^{\ast }d\left( x_{n},p\right) +\alpha _{n}\varphi
_{2n}.
\end{split}
\label{2.2}
\end{equation}%
Substituting (2.2) in (2.1), gives that

\begin{eqnarray*}
d\left( x_{n+1},p\right) &\leq &d\left( x_{n},p\right) +\alpha
_{n}k_{2n}M_{2}^{\ast }d\left( x_{n},p\right) +\alpha _{n}k_{2n}k_{1n\
}M_{2}^{\ast }M_{1}^{\ast }d\left( x_{n},p\right) \\
&&+\alpha _{n}k_{1n}M_{1}^{\ast }d\left( x_{n},p\right) +\alpha
_{n}k_{2n}\xi _{2}\left( M_{2}\right) +\alpha _{n}k_{1n}\xi _{1}\left(
M_{1}\right) \\
&&+\alpha _{n}k_{2n}k_{1n\ }M_{1}^{\ast }\xi _{2}\left( M_{2}\right) +\alpha
_{n}k_{1n}\varphi _{2n}M_{1}^{\ast } \\
&&+\alpha _{n}\varphi _{2n}+\alpha _{n}\varphi _{1n}.
\end{eqnarray*}%
\newline
Thus, we have%
\begin{equation*}
d\left( x_{n+1},p\right) \leq \left( 1+\left( k_{1n\ }+k_{2n\ }\right)
a_{2}\right) d\left( x_{n},p\right) +\left( k_{2n\ }+k_{1n\ }+\varphi
_{2n}+\varphi _{1n}\right) a_{2}
\end{equation*}%
for some constant $a_{2}>0.$\newline
Continuing in a similar fashion, we obtain that%
\begin{equation*}
d\left( x_{n+1},p\right) \leq \left( 1+a\sum\limits_{j=1}^{m}k_{jn\ }\right)
d\left( x_{n},p\right) +a\sum\limits_{j=1}^{m}\left( k_{jn\ }+\varphi
_{jn}\right)
\end{equation*}%
for some constant $a>0.$\newline
Appealing to Lemma 2.4, the above estimate implies that the sequence $%
\left\{ x_{n}\right\} $ is bounded and $\lim\limits_{n\rightarrow \infty
}d\left( x_{n},p\right) $ exists. This completes the proof.\newline
\textbf{Lemma 2.8.} \textit{Let }$K$ \textit{be a nonempty closed and convex
subset of a uniformly convex hyperbolic space }$X$\textit{\ with monotone
modulus of uniform convexity }$\eta $ \textit{and let }$\left\{
T_{i}\right\} _{i=1}^{m}:K\rightarrow K$ \textit{be a finite family of
uniformly continuous total asymptotically nonexpansive mappings with
sequences }$\{k_{in}\}$ and $\{\varphi _{in}\},$ $n\geq 1,\ i=1,2,\cdots ,m$ 
\textit{such that }$F:=\cap _{i=1}^{m}F\left( T_{i}\right) \neq \emptyset .$ 
\textit{For} $i=1,2,\cdots ,m$ \textit{if the following conditions are
satisfied:}\newline
(C1) $\sum\limits_{n=1}^{\infty }k_{in}<\infty $ \textit{and}$\
\sum\limits_{n=1}^{\infty }\varphi _{in}<\infty ;$\newline
(C2) \textit{there exists constants }$M_{i},\ M_{i}^{\ast }>0$ \textit{such
that }$\zeta _{i}\left( \lambda _{i}\right) \leq M_{i}^{\ast }\lambda _{i}$%
\textit{\ for all }$\lambda _{i}\geq M_{i}$\textit{. Then }$%
\lim\limits_{n\rightarrow \infty }d\left( x_{n},T_{i}x_{n}\right) =0,$ $%
i=1,2,\ldots ,m.$\newline
\textbf{Proof.} It follows from Lemma 2.7 that $\lim\limits_{n\rightarrow
\infty }d\left( x_{n},p\right) $ exists. Without loss of generality, we can
assume that $\lim\limits_{n\rightarrow \infty }d\left( x_{n},p\right) =r>0.$
We first distinguish two cases to show that $\lim\limits_{n\rightarrow
\infty }d\left( x_{n},T_{i}^{n}\right) =0,$ $i=1,2,\ldots ,m.$\newline
\textbf{Case 1. }For $m=1$, we proceed as follows:%
\begin{equation}
\lim\limits_{n\rightarrow \infty }d\left( x_{n+1},p\right)
=\lim\limits_{n\rightarrow \infty }d\left( W\left(
x_{n},T_{1}^{n}x_{n},\alpha _{n}\right) ,p\right) =r.  \label{2.3}
\end{equation}%
It follows from the definition that%
\begin{eqnarray*}
d\left( T_{1}^{n}x_{n},p\right) &=&d\left( T_{1}^{n}x_{n},T_{1}^{n}p\right)
\leq d\left( x_{n},p\right) +k_{1n\ }\zeta _{1}(d\left( x_{n},p\right)
)+\varphi _{1n} \\
&\leq &d\left( x_{n},p\right) +k_{1n\ }M_{1}^{\ast }d\left( x_{n},p\right)
+\alpha _{n}k_{1n}\zeta _{1}\left( M_{1}\right) +\varphi _{1n} \\
&=&(1+k_{1n\ }M_{1}^{\ast })d\left( x_{n},p\right) +\alpha _{n}k_{1n}\zeta
_{1}\left( M_{1}\right) +\varphi _{1n}.
\end{eqnarray*}%
Taking $\lim \sup $ on both sides of the above estimate and utilizing the
fact that $k_{1n},\varphi _{1n}\overset{n\rightarrow \infty }{\rightarrow }%
0, $ we get%
\begin{equation*}
\limsup\limits_{n\rightarrow \infty }d\left( T_{1}^{n}x_{n},p\right) \leq r.
\end{equation*}%
Moreover, $\limsup\limits_{n\rightarrow \infty }d\left( x_{n},p\right) \leq
r.$ Hence, the conclusion, i.e., $\lim\limits_{n\rightarrow \infty }d\left(
T_{1}^{n}x_{n},x_{n}\right) =0$ follows from Lemma 2.5.\newline
\textbf{Case 2.} For $m=2,$ iteration (1.1) reduces to:%
\begin{eqnarray*}
x_{n+1} &=&W\left( x_{n},T_{1}^{n}y_{1n},\alpha _{n}\right) \\
y_{1n} &=&W\left( x_{n},T_{2}^{n}x_{n},\alpha _{n}\right) .
\end{eqnarray*}%
We first calculate the following estimates:%
\begin{eqnarray*}
d\left( y_{1n},p\right) &=&d\left( W\left( x_{n},T_{2}^{n}x_{n},\alpha
_{n}\right) ,p\right) \\
&\leq &\left( 1-\alpha _{n}\right) d\left( x_{n},p\right) +\alpha
_{n}d\left( T_{2}^{n}x_{n},p\right) \\
&\leq &\left( 1-\alpha _{n}\right) d\left( x_{n},p\right) +\alpha _{n}\left[
d\left( x_{n},p\right) +k_{2n\ }\zeta _{2}\left( d\left( x_{n},p\right)
\right) +\varphi _{2n}\right] \\
&=&d\left( x_{n},p\right) +\alpha _{n}k_{2n\ }M_{2}^{\ast }d\left(
x_{n},p\right) +\alpha _{n}k_{2n}\zeta _{2}\left( M_{2}\right) +\alpha
_{n}\varphi _{2n} \\
&=&(1+k_{2n\ }M_{2}^{\ast })d\left( x_{n},p\right) +\alpha _{n}k_{2n}\zeta
_{2}\left( M_{2}\right) +\alpha _{n}\varphi _{2n}.
\end{eqnarray*}%
Taking $\lim \sup $ on both sides of the above estimate, we have%
\begin{equation}
\limsup\limits_{n\rightarrow \infty }d\left( y_{1n},p\right) \leq r.
\label{2.4}
\end{equation}%
Now observe that%
\begin{equation*}
d\left( T_{1}^{n}y_{1n},p\right) \leq (1+k_{1n\ }M_{1}^{\ast })d\left(
y_{1n},p\right) +k_{1n}\zeta _{1}\left( M_{1}\right) +\varphi _{1n}.
\end{equation*}%
Taking $\lim \sup $ on both sides of the above estimate\ and utilizing
(2.4), we get%
\begin{equation}
\limsup\limits_{n\rightarrow \infty }d\left( T_{1}^{n}y_{1n},p\right) \leq r.
\label{2.5}
\end{equation}%
Since $\limsup\limits_{n\rightarrow \infty }d\left( x_{n},p\right) \leq r,$
therefore now it follows from (2.5) and Lemma 2.5, that%
\begin{equation}
\lim\limits_{n\rightarrow \infty }d\left( x_{n},T_{1}^{n}y_{1n}\right) =0.
\label{2.6}
\end{equation}%
Observe that%
\begin{eqnarray*}
d\left( x_{n},p\right) &\leq &d\left( x_{n},T_{1}^{n}y_{1n}\right) +d\left(
T_{1}^{n}y,p\right) \\
&\leq &d\left( x_{n},T_{1}^{n}y_{1n}\right) +(1+k_{1n\ }M_{1}^{\ast
})d\left( y_{1n},p\right) +k_{1n}\zeta _{1}\left( M_{1}\right) +\varphi
_{1n}.
\end{eqnarray*}%
Hence, we deduce from the above estimate that%
\begin{equation}
r\leq \liminf\limits_{n\rightarrow \infty }d\left( y_{1n},p\right) .
\label{2.7}
\end{equation}%
The estimates (2.4) and (2.7) collectively imply that 
\begin{equation}
\lim\limits_{n\rightarrow \infty }d\left( y_{1n},p\right) =d\left( W\left(
x_{n},T_{2}^{n}x_{n},\alpha _{n}\right) ,p\right) =r.  \label{2.8}
\end{equation}%
Utilizing the total asymptotically nonexpansiveness of $T_{2}$ and the fact
that $k_{2n},\varphi _{2n}\overset{n\rightarrow \infty }{\rightarrow }0,$ we
have%
\begin{equation}
\limsup_{n\rightarrow \infty }d\left( T_{2}^{n}x_{n},p\right) \leq r.
\label{2.9}
\end{equation}%
Now (2.8), (2.9) and Lemma 2.5, imply that%
\begin{equation}
\lim\limits_{n\rightarrow \infty }d\left( T_{2}^{n}x_{n},x_{n}\right) =0.
\label{2.10}
\end{equation}%
Note that $d\left( y_{1n},x_{n}\right) \leq \Theta \cdot d\left(
T_{2}^{n}x_{n},x_{n}\right) ,$ for some $\Theta >0.$ This implies - - on
utilizing (2.10) - - that%
\begin{equation}
\lim\limits_{n\rightarrow \infty }d\left( y_{1n},x_{n}\right) =0.
\label{2.11}
\end{equation}%
It follows from (2.6) and (2.11) that%
\begin{equation}
\lim\limits_{n\rightarrow \infty }d\left( y_{1n},T_{1}^{n}y_{1n}\right) =0.
\label{2.12}
\end{equation}%
Moreover%
\begin{equation*}
d\left( T_{1}^{n}x_{n},x_{n}\right) \leq d\left(
T_{1}^{n}x_{n},T_{1}^{n}y_{1n}\right) +d\left( T_{1}^{n}y_{1n},x_{n}\right) .
\end{equation*}%
Since $T_{1}$ is uniformly continuous, therefore letting $n\rightarrow
\infty $ in the above estimate and utilizing (2.6) and (2.11), we get%
\begin{equation*}
\lim\limits_{n\rightarrow \infty }d\left( T_{1}^{n}x_{n},x_{n}\right) =0.
\end{equation*}%
Hence%
\begin{equation*}
\lim\limits_{n\rightarrow \infty }d\left( T_{i}^{n}x_{n},x_{n}\right) =0%
\text{ for }i=1,2.
\end{equation*}%
Continuing in a similar way for $i=1,2,\ldots ,m,$ we get that%
\begin{equation*}
\lim\limits_{n\rightarrow \infty }d\left( T_{i}^{n}y_{in},x_{n}\right)
=\lim\limits_{n\rightarrow \infty }d\left( y_{in},x_{n}\right)
=\lim\limits_{n\rightarrow \infty }d\left( T_{i}^{n}y_{in},y_{in}\right)
=\lim\limits_{n\rightarrow \infty }d\left( T_{i}^{n}x_{n},x_{n}\right) =0.
\end{equation*}%
Now, utilizing the uniform contiuity of $T_{i},$ the following estimate:%
\begin{equation*}
d\left( x_{n},T_{i}x_{n}\right) \leq d\left( x_{n},T_{i}^{n}x_{n}\right)
+d\left( T_{i}^{n}x_{n},T_{i}x_{n}\right)
\end{equation*}%
implies that%
\begin{equation*}
\lim\limits_{n\rightarrow \infty }d\left( T_{i}x_{n},x_{n}\right) =0\text{
for }i=1,2,\ldots ,m.
\end{equation*}%
This completes the proof.

\section{Existence and Convergence Results}

In this section, we first establish the existence of fixed point of total
asymptotically nonexpansive mappings. Our proof closely follows \cite[%
Theorem 3.3]{Kohlenbach JEMS} for the existence of fixed point of
asymptotically nonexpansive mappings in uniformly convex hyperbolic spaces.%
\newline
\textbf{Theorem 3.1.} \textit{Let }$(X,d,W)$\textit{\ be a complete
uniformly convex hyperbolic space having a monotone modulus of uniform
convexity }$\eta $\textit{. Let }$K$\textit{\ be a nonempty closed bounded
and convex subset of }$X.$ \textit{Then every mapping of total
asymptotically nonexpansive }$T:K\rightarrow K$\textit{\ has a fixed point.}%
\newline
\textbf{Proof.} For any $y\in K$, define%
\begin{equation*}
K_{y}=\left\{ \tau \in 
\mathbb{R}
^{+}:\text{ there exists }x\in K,n_{0}\in 
\mathbb{N}
\text{ such that }d(T^{i}y,x)\leq \tau \text{ for all }i\geq n_{0}\right\} .
\end{equation*}%
\newline
It is easy to see that $K_{y}$ is nonempty as $diam(K)\in K_{y}.$ Let $%
\alpha _{y}:=\inf K_{y}$ then for any $\theta >0$, there exists $\tau
_{\theta }\in K_{y}$ such that $\tau _{\theta }<\alpha _{y}+\theta ,$ so
there exists $x\in K$ and $n_{0}\in 
\mathbb{N}
$ such that%
\begin{equation}
d(T^{i}y,x)\leq \tau _{\theta }<\alpha _{y}+\theta \text{ for all }i\geq
n_{0}.  \label{3.1}
\end{equation}%
Obviously, $\alpha _{y}\geq 0.$ Here, we distinguish two cases:\newline
\textbf{Case 1.} $\alpha _{y}=0.$ We show that the sequence $\left\{
T^{i}y\right\} _{i=1}^{\infty }$ is Cauchy. Let $m,n\geq n_{0}\in 
\mathbb{N}
$ and let $\epsilon >0,$ then applying (3.1) with $\theta =\frac{\epsilon }{2%
},$ we get%
\begin{equation*}
d(T^{m}y,T^{n}y)\leq d(T^{m}y,x)+d(T^{n}y,x)<\frac{\epsilon }{2}+\frac{%
\epsilon }{2}=\epsilon .
\end{equation*}%
Hence $\left\{ T^{i}y\right\} _{i=1}^{\infty }$ is Cauchy and consequently
converges to some $z\in K.$ Using the Definition of $T$ choose $n_{1}\in 
\mathbb{N}
$ such that%
\begin{equation*}
\left\{ \left( d(T^{i}z,T^{i}x)-d(z,x)-k_{n}\xi \left( d(z,x)\right)
-\varphi _{n}\right) \right\} \leq \frac{\epsilon }{2}\text{ for all }i\geq
n_{1}.
\end{equation*}%
Since $T^{i}y\rightarrow z,$ this implies $d(T^{m}y,z),d(T^{m-i}y,z)\leq 
\frac{\epsilon }{2}$ for all $m>i.$ Thus, if $i\geq m$ then%
\begin{eqnarray*}
d(z,T^{i}z) &\leq &d(z,T^{m}y)+d(T^{m}y,T^{i}z) \\
&\leq &d(z,T^{m}y)+\left( d(T^{i}z,T^{i}\left( T^{m-i}y\right) )-k_{i}\xi
\left( d(T^{m-i}y,z)\right) -\varphi _{i}\right) \\
&&+k_{n}\xi \left( d(T^{m-i}y,z)\right) +\varphi _{n} \\
&\leq &d(z,T^{m}y)+\left( d(T^{i}z,T^{i}\left( T^{m-i}y\right) )-k_{i}\xi
\left( d(T^{m-i}y,z)\right) -\varphi _{i}\right) \\
&&+k_{i}M^{\ast }d(T^{m-i}y,z)+\varphi _{i} \\
&\leq &\epsilon +k_{i}M^{\ast }\frac{\epsilon }{2}+\varphi _{i},
\end{eqnarray*}%
where $M^{\ast }>0$ such that $\zeta \left( \lambda \right) \leq M^{\ast
}\lambda $.\newline
Letting $i\rightarrow \infty $ in the above estimate, we have $%
T^{i}z\rightarrow z.$ Hence the continuity of $T$ implies that $z$ is a
fixed point of $T,i.e.,$%
\begin{equation*}
Tz=T(T^{i}z)=T^{i+1}z=z\text{ as }i\rightarrow \infty
\end{equation*}%
\newline
\textbf{Case 2.} $\alpha _{y}>0.$ For any $n\geq 1$, we define%
\begin{equation*}
C_{n}=\tbigcup_{k\geq 1}\tbigcap_{i\geq k}\overline{U}\left( T^{i}y,\alpha
_{y}+\frac{1}{n}\right) ,~~~D_{n}:=\overline{C_{n}}\cap C_{n}.
\end{equation*}%
By (3.1) with $\theta =\frac{1}{n},$ there exist $x\in K$ and $k\geq 1$ such
that $x\in \tbigcap_{i\geq k}\overline{U}\left( T^{i}y,\alpha _{y}+\frac{1}{n%
}\right) ;$ hence $D_{n}$ is nonempty. Moreover, $\{D_{n}\}$ is a decreasing
sequence of nonempty-bounded closed convex subsets of $X$, hence, we can
apply Proposition 2.2 to derive that%
\begin{equation*}
D:=\tbigcap_{n\geq 1}D_{n}\neq \emptyset .
\end{equation*}%
For any $x\in D$ and $\theta >0,$ take $n^{\ast }\in 
\mathbb{N}
$ such that $\frac{2}{n^{\ast }}\leq \theta .$ Since $x\in \overline{%
C_{n^{\ast }}},$ so there exists a sequence $\{x_{(n^{\ast })n}\}$ in $%
C_{n^{\ast }}$ such that $\lim_{n\rightarrow \infty }x_{(n^{\ast })n}=x.$
Let $n^{^{^{\prime }}}\geq 1$ be such that $d\left( x_{(n^{\ast
})n},x\right) \leq \frac{1}{n^{\ast }}$ for all $n\geq n^{^{\prime }}$ and $%
s\geq 1$ such that $x_{(n^{\ast })n^{^{\prime }}}\in \tbigcap_{i\geq s}%
\overline{U}\left( T^{i}y,\alpha _{y}+\frac{1}{n^{\ast }}\right) .$ It
follows that for all $i\geq s$%
\begin{equation}
d(T^{i}y,x)\leq d\left( T^{i}y,x_{(n^{\ast })n^{^{\prime }}}\right) +d\left(
x_{(n^{\ast })n^{^{\prime }}},x\right) \leq \alpha _{y}+\frac{1}{n^{\ast }}+%
\frac{1}{n^{\ast }}\leq \alpha _{y}+\theta .  \label{3.2}
\end{equation}%
In the sequel, we shall prove that any point of $D$ is a fixed point of $T$.
We assume towards contradiction, that is $Tx\neq x$ for any $x\in D.$ Then $%
T^{n}x\nrightarrow x($by Case 1$),$so there exists $\epsilon >0$ and $%
\widehat{n}\in 
\mathbb{N}
$ such that%
\begin{equation}
d\left( T^{n}x,x\right) \geq \epsilon \text{ for all }n\geq \widehat{n}.
\label{3.3}
\end{equation}%
We can assume that $\epsilon \in (0,2]$ so that $\frac{\epsilon }{\alpha
_{y}+1}\in (0,2].$ Hence there exists $\theta _{y}\in (0,1]$ such that%
\begin{equation*}
1-\eta \left( \alpha _{y}+1,\frac{\epsilon }{\alpha _{y}+1}\right) \leq 
\frac{\alpha _{y}-\theta _{y}}{\alpha _{y}+\theta _{y}}.
\end{equation*}%
Observe that $\lim_{n\rightarrow \infty }\left( \left( 1+k_{n}M^{\ast
}\right) \left( \alpha _{y}+\frac{\theta _{y}}{2}\right) +\varphi
_{n}\right) =\alpha _{y}+\frac{\theta _{y}}{2}<\alpha _{y}+\theta _{y},$ so
there exists $\widetilde{n}\in 
\mathbb{N}
$ such that%
\begin{equation}
\left( 1+k_{n}M^{\ast }\right) \left( \alpha _{y}+\frac{\theta _{y}}{2}%
\right) +\varphi _{n}<\alpha _{y}+\theta _{y}\text{ for all }n\geq 
\widetilde{n}.  \label{3.4}
\end{equation}%
\newline
Applying (3.2) with $\theta =\frac{\theta _{y}}{2}$, there exists $n^{\ast
\ast }\in 
\mathbb{N}
$ such that%
\begin{equation}
d\left( T^{i}y,x\right) \leq \alpha _{y}+\frac{\theta _{y}}{2},\text{ for
all }i\geq n^{\ast \ast }.  \label{3.5}
\end{equation}%
By the Definition of $T$ choose $n^{^{^{\prime \prime }}}\in 
\mathbb{N}
$ such that%
\begin{equation*}
\left\{ \left( d(T^{m}z,T^{m}x)-d(z,x)-k_{n}\xi \left( d(z,x)\right)
-\varphi _{n}\right) \right\} \leq \frac{\theta _{y}}{2}\text{ for all }%
m\geq n^{^{\prime \prime }}.
\end{equation*}%
Applying (3.3) with $\widehat{n}:=\widetilde{n},$ we get%
\begin{equation}
d\left( T^{n^{\ast }}x,x\right) \geq \epsilon \geq \frac{\epsilon }{\alpha
_{y}+1}\left( \alpha _{y}+\theta _{y}\right) \text{ for all }n^{\ast }\geq 
\widetilde{n}.  \label{3.6}
\end{equation}%
Let now $m\in 
\mathbb{N}
$ be such that $m\geq n^{\ast }+n^{\ast \ast },$ then the estimate (3.5)
implies that%
\begin{equation}
d\left( T^{m}y,x\right) \leq \alpha _{y}+\frac{\theta _{y}}{2}<\alpha
_{y}+\theta _{y}.  \label{3.7}
\end{equation}%
Moreover, observe that%
\begin{eqnarray*}
d\left( T^{n^{\ast }}x,T^{m}y\right) &=&d\left( T^{n^{\ast }}x,T^{n^{\ast
}}(T^{m-n^{\ast }}y)\right) \\
&\leq &d\left( x,T^{m-n^{\ast }}y\right) +k_{n^{\ast }}\xi \left(
d(x,T^{m-n^{\ast }}y)\right) +\varphi _{n^{\ast }} \\
&\leq &d\left( x,T^{m-n^{\ast }}y\right) +k_{n^{\ast }}M^{\ast }d\left(
x,T^{m-n^{\ast }}y\right) +\varphi _{n^{\ast }} \\
&\leq &\left( 1+k_{n^{\ast }}M^{\ast }\right) \left( \alpha _{y}+\frac{%
\theta _{y}}{2}\right) +\varphi _{n^{\ast }}.
\end{eqnarray*}%
Utilizing (3.4) in the above estimate, we get%
\begin{equation}
d\left( T^{n^{\ast }}x,T^{m}y\right) <\alpha _{y}+\theta _{y}.  \label{3.8}
\end{equation}%
It follows from the estimates (3.6)-(3.8) and uniform convexity of $X$, that%
\begin{eqnarray*}
d\left( W\left( x,T^{n^{\ast }}x,\frac{1}{2}\right) ,T^{m}y\right) &\leq
&\left( 1-\eta \left( \alpha _{y}+\theta _{y},\frac{\epsilon }{\alpha _{y}+1}%
\right) \right) \left( \alpha _{y}+\theta _{y}\right) \\
&\leq &\left( 1-\eta \left( \alpha _{y}+1,\frac{\epsilon }{\alpha _{y}+1}%
\right) \right) \left( \alpha _{y}+\theta _{y}\right) \text{ \ \ (}\eta 
\text{ is monotone)} \\
&\leq &\frac{\alpha _{y}-\theta _{y}}{\alpha _{y}+\theta _{y}}\cdot \left(
\alpha _{y}+\theta _{y}\right) =\alpha _{y}-\theta _{y}.
\end{eqnarray*}%
Hence, there exists $\widehat{n}:=n^{\ast }+n^{\ast \ast }$ and $z:=W\left(
x,T^{n^{\ast }}x,\frac{1}{2}\right) \in K$ such that%
\begin{equation*}
d\left( z,T^{m}y\right) \leq \alpha _{y}-\theta _{y}\text{ for all }m\geq 
\widehat{n}.
\end{equation*}%
This implies that $\alpha _{y}-\theta _{y}\in K_{y}$ which contradict the
fact that $\alpha _{y}=\inf ~K_{y}.$ Hence $x$ is a fixed point of $T.$ This
completes the proof.

We now use the concept of asymptotic center of a bounded sequence to
strengthen the above existence result.\newline
\textbf{Theorem 3.2.} \textit{Let }$(X,d,W)$\textit{\ be a complete
uniformly convex hyperbolic space having a monotone modulus of uniform
convexity }$\eta $\textit{. Let }$K$\textit{\ be a nonempty closed bounded
and convex subset of }$X$\textit{\ and let }$T:K\rightarrow K$\textit{\ be a
total asymptotically nonexpansive mapping. If }$\{T^{n}x\}$\textit{\ is
bounded for some }$x\in K$\textit{\ and }$z\in A_{K}(\{T^{n}x\})$\textit{\
then }$z$\textit{\ is a fixed point of }$T.$\textit{\newline
}\textbf{Proof.} Let $\{x_{n}=T^{n}x\}$ is bounded and let $z\in
A_{K}(\{x_{n}\})$. It follows from Lemma 2.1 that $z$ is the unique
asymptotic center of $\{x_{n}\}$. Assume $\{y_{m}=T^{m}z\}\in K.$ For
integer $n>m\geq 1,$ we have%
\begin{eqnarray*}
d\left( x_{n},y_{m}\right)  &=&d\left( T^{m}x_{n-m},T^{m}z\right)  \\
&\leq &d\left( x_{n-m},z\right) +k_{m}\xi \left( d\left( x_{n-m},z\right)
\right) +\varphi _{m} \\
&\leq &\left( 1+k_{m}M^{\ast }\right) d\left( x_{n-m},z\right) +k_{m}\xi
\left( M\right) +\varphi _{m}.
\end{eqnarray*}%
Letting $m\rightarrow \infty $ in the above estimate, we get%
\begin{eqnarray*}
r\left( x_{n},\{y_{m}\}\right)  &=&\limsup_{m\rightarrow \infty }d\left(
x_{n},y_{m}\right)  \\
&\leq &\limsup_{m\rightarrow \infty }d\left( x_{n-m},z\right)  \\
&=&A_{K}(\{x_{n}\}).
\end{eqnarray*}%
Now, Lemma 2.6 implies that $y_{m}\rightarrow z$ as $m\rightarrow \infty .$
The continuity of $T$ implies that $z$ is a fixed point of $T.$ This
completes the proof.\medskip 

Rest of the paper deals with the convergence analysis of iteration (1.1) for
the approximation of common fixed points of a finite family of total
asymptotically nonexpansive mappings.\newline
\textbf{Theorem 3.3}.\textbf{\ }\textit{Let }$K$ \textit{be a nonempty
closed and convex subset of a uniformly convex hyperbolic space }$X$\textit{%
\ with monotone modulus of uniform convexity }$\eta $ \textit{and let }$%
\left\{ T_{i}\right\} _{i=1}^{m}:K\rightarrow K$ \textit{be a finite family
of uniformly continuous total asymptotically nonexpansive mappings with
sequences }$\{k_{in}\}$\textit{\ and} $\{\varphi _{in}\},$ $n\geq 1,\
i=1,2,\cdots ,m$ \textit{such that }$F:=\cap _{i=1}^{m}F\left( T_{i}\right)
\neq \emptyset .$ \textit{For} $i=1,2,\cdots ,m$ \textit{if the following
conditions are satisfied:}\newline
(C1) $\sum\limits_{n=1}^{\infty }k_{in}<\infty $ \textit{and}$\
\sum\limits_{n=1}^{\infty }\varphi _{in}<\infty ;$\newline
(C2) \textit{there exists constants }$M_{i},\ M_{i}^{\ast }>0$ \textit{such
that }$\zeta _{i}\left( \lambda _{i}\right) \leq M_{i}^{\ast }\lambda _{i}$%
\textit{\ for all }$\lambda _{i}\geq M_{i}$\textit{. Then the sequence }$%
\{x_{n}\}$\textit{\ defined in (1.1) }$\triangle $-\textit{converges to a
common fixed point in }$F.\smallskip $\textit{\newline
}\textbf{Proof.} It follows from Lemma 2.7 that $\{x_{n}\}$ is bounded.
Therefore by Lemma 2.1, $\{x_{n}\}$ has a unique asymptotic center, that is, 
$A_{K}(\{x_{n}\})=\{x\}.$ Let $\{u_{n}\}$ be any subsequence of $\{x_{n}\}$
such that $A_{K}(\{u_{n}\})=\{u\}$ and%
\begin{equation}
\lim_{n\rightarrow \infty }d(u_{n},T_{i}u_{n})=0\text{ for all }i=1,2,\cdots
,m.  \label{3.9}
\end{equation}%
Next, we show that $u\in F.$ For each $j\in \{1,2,3,\cdots ,m\},$ we define
a sequence $\{z_{n}\}$ in $K$ by $z_{i}=T_{j}^{i}u.$ So we calculate%
\begin{eqnarray*}
d(z_{i},u_{n}) &\leq
&d(T_{j}^{i}u,T_{j}^{i}u_{n})+d(T_{j}^{i}u_{n},T_{j}^{i-1}u_{n})+\cdots
+d(T_{j}u_{n},u_{n}) \\
&\leq &\left( 1+k_{jn}M_{j}^{\ast }\right) d(u,u_{n})+k_{jn}\xi _{j}\left(
M_{j}\right) +\sum_{i=0}^{r-1}d(T_{j}^{i}u_{n},T_{j}^{i+1}u_{n}).
\end{eqnarray*}%
Taking $\lim \sup $ on both sides of the above estimate and utilizing (3.9)
and the fact that each $T_{j}$ is uniformly continuous, we have%
\begin{equation*}
r(z_{i},\{u_{n}\})=\limsup_{n\rightarrow \infty }d(z_{i},u_{n})\leq
\limsup_{n\rightarrow \infty }d(u,u_{n})=r(u,\{u_{n}\}).
\end{equation*}%
This implies that $\left\vert r(z_{i},\{u_{n}\})-r(u,\{u_{n}\})\right\vert
\rightarrow 0$ as $i\rightarrow \infty .$ It follows from Lemma 2.6 that $%
\lim_{i\rightarrow \infty }T_{j}^{i}u=u.$ Again, utilizing the uniform
continuity of $T_{j},$ we have that $T_{j}(u)=T_{j}(\lim_{i\rightarrow
\infty }T_{j}^{i}u)=\lim_{i\rightarrow \infty }T_{j}^{i+1}u=u.$ From the
arbitrariness of $j,$ we conclude that $u$ is the common fixed point of $%
\{T_{j}\}_{j=1}^{m}.$ It remains to show that $x=u.$ In fact, uniqueness of
asymptotic center implies that%
\begin{eqnarray*}
\limsup_{n\rightarrow \infty }d(u_{n},u) &<&\limsup_{n\rightarrow \infty
}d(u_{n},x) \\
&\leq &\limsup_{n\rightarrow \infty }d(x_{n},x) \\
&<&\limsup_{n\rightarrow \infty }d(x_{n},u) \\
&=&\limsup_{n\rightarrow \infty }d(u_{n},u).
\end{eqnarray*}%
This is a contradiction. Hence $x=u.$ This implies that $u$ is the unique
asymptotic center of $\{x_{n}\}$ for every subsequence $\{u_{n}\}$ of $%
\{x_{n}\}.$ This completes the proof.\medskip 

The strong convergence of iteration (1.1) can easily be established under
compactness condition of $K$ or $T(K).$ As a consequence, we can get a
generalized version of \cite[Theorem 12]{Chidume JMAA} and \cite[Theorem 3.3]%
{KDF JMAA} to the general setup of uniformly convex hyperbolic spaces,
respectively. Next, we give a necessary and sufficient condition for the
strong convergence of iteration (1.1).\newline
\textbf{Theorem 3.4}.\textbf{\ }\textit{Let }$K$ \textit{be a nonempty
closed and convex subset of a uniformly convex hyperbolic space }$X$\textit{%
\ with monotone modulus of uniform convexity }$\eta $ \textit{and let }$%
\left\{ T_{i}\right\} _{i=1}^{m}:K\rightarrow K$ \textit{be a finite family
of uniformly continuous total asymptotically nonexpansive mappings with
sequences }$\{k_{in}\}$\textit{\ and} $\{\varphi _{in}\},$ $n\geq 1,\
i=1,2,\cdots ,m.$ \textit{Suppose that }$F:=\cap _{i=1}^{m}F\left(
T_{i}\right) \neq \emptyset $\textit{\ and the following conditions are
satisfied:}\newline
(C1) $\sum\limits_{n=1}^{\infty }k_{in}<\infty $ \textit{and}$\
\sum\limits_{n=1}^{\infty }\varphi _{in}<\infty ;$\newline
(C2) \textit{there exists constants }$M_{i},\ M_{i}^{\ast }>0$ \textit{such
that }$\zeta _{i}\left( \lambda _{i}\right) \leq M_{i}^{\ast }\lambda _{i}$%
\textit{\ for all }$\lambda _{i}\geq M_{i}$\textit{. Then the sequence }$%
\{x_{n}\}$\textit{\ defined in (1.1) converges strongly to a common fixed
point in }$F$\textit{\ if and only if }$\lim \inf_{n\rightarrow \infty
}d(x_{n},F)=0.$\newline
\textbf{Proof.} The necessity of the conditions is obvious. Thus, we only
prove the sufficiency. It follows from Lemma 2.8 that $\left\{
d(x_{n},p)\right\} _{n=1}^{\infty }$\ converges. Moreover, $\lim
\inf_{n\rightarrow \infty }d(x_{n},F)=0$ implies that $\lim_{n\rightarrow
\infty }d(x_{n},F)=0.$ This completes the proof.\medskip

It is remarked that there are certain situations when the domain $D(T)$ of a
nonlinear mapping $T$ is a proper subset of the under laying space $X.$ In
such situations, the iterative schema for the approximation of fixed points
of $T$ may failed to be well-defined. It is therefore natural to study
non-self behaviour of the nonlinear mappings.

We recall that a non-self-mapping $T:K\rightarrow X$ is called total
asymptotically nonexpansive mappings if there exists nonnegative real
sequences $\{k_{n}\}$ and $\{\varphi _{n}\}$ with $k_{n}\overset{%
n\rightarrow \infty }{\rightarrow }0,\varphi _{n}\overset{n\rightarrow
\infty }{\rightarrow }0$ and a strictly increasing continuous function $\xi :%
\mathbb{R}
^{+}\rightarrow 
\mathbb{R}
^{+}$ with $\xi (0)=0,$ such that 
\begin{equation*}
d(T\left( PT\right) ^{n-1}x,T\left( PT\right) ^{n-1}y)\leq d(x,y)+k_{n}\xi
\left( d(x,y)\right) +\varphi _{n}\text{ \ for all \ }x,y\in K,~n\geq 1,
\end{equation*}%
where $P:X\rightarrow K$ is a nonexpansive retract.

Hence, one can establish the strong and $\triangle $-convergence results --
as in Theorems 3.3-3.4 with a slight modification -- for the following
iteration schema involving a finite family of totaly asymptotically
nonexpansive non-self mappings:%
\begin{equation*}
\begin{array}{l}
x_{1}\in K, \\ 
x_{n+1}=P\left( W\left( x_{n},T_{1}\left( PT_{1}\right) ^{n-1}x_{n},\alpha
_{n}\right) \right) \text{, if }m=1,n\geq 1, \\ 
x_{1}\in K, \\ 
x_{n+1}=P\left( W\left( x_{n},T_{1}\left( PT_{1}\right) ^{n-1}y_{1n},\alpha
_{n}\right) \right) , \\ 
y_{1n}=P\left( W\left( x_{n},T_{2}\left( PT_{2}\right) ^{n-1}y_{2n},\alpha
_{n}\right) \right) , \\ 
y_{2n}=P\left( W\left( x_{n},T_{3}\left( PT_{3}\right) ^{n-1}y_{3n},\alpha
_{n}\right) \right) , \\ 
\multicolumn{1}{c}{\vdots} \\ 
y_{(m-2)n}=P\left( W\left( x_{n},T_{m-1}\left( PT_{m-1}\right)
^{n-1}y_{\left( m-1\right) n},\alpha _{n}\right) \right) , \\ 
y_{(m-1)n}=P\left( W\left( x_{n},T_{m}\left( PT_{m}\right)
^{n-1}y_{mn},\alpha _{n}\right) \right) \text{, if }m\geq 2,n\geq 1.%
\end{array}%
\end{equation*}%
\textbf{Remark 3.5.} (i) It is worth to mention that Theorems 3.3-3.4 can
easily be extended to the class of mappings with bounded error terms;\newline
(ii) Lemma 2.8 improves and generalizes \cite[Theorem 7]{Chidume JMAA} and 
\cite[Lemma 3.1(i)]{KDF JMAA} to the general setup of uniformly convex
hyperbolic spaces respectively;\newline
(iii) Lemma 2.9 improves and generalizes \cite[Theorem 11]{Chidume JMAA} and 
\cite[Lemma 3.1(iii)]{KDF JMAA} to the setting as mentioned in (ii);\newline
(iv) Theorems 3.3 improves \cite[Theorem 3.5]{Chang AMC} and \cite[Theorem
3.1]{Khan JIA} for a finite family of total asymptotically nonexpansive
mappings to the setup of spaces as mentioned in (ii);\newline
(v) Theorems 3.3 sets analogue of \cite[Theorem 3.2]{KDF JMAA} to the
setting as mentioned in (ii);\newline
(vi) Theorem 3.4 generalizes \cite[Theorem 8]{Chidume JMAA} and \cite[%
Theorem 2.2]{KDF JMAA} to the setting as defined in (ii).

\end{document}